\documentstyle[12pt]{article}

\newlength{\abstractwidth}
\renewcommand{\thefootnote}{\fnsymbol{footnote}}
\renewcommand{\thanks}[1]{\footnote{#1}} 
\newcommand{\starttext}{ \setcounter{footnote}{0}
\renewcommand{\thefootnote}{\arabic{footnote}}}

\newcommand{\be}{\begin{equation}}
\newcommand{\bea}{\begin{eqnarray}}
\newcommand{\eea}{\end{eqnarray}} \newcommand{\ee}{\end{equation}}
 
\renewcommand{\>}{\rangle}
\def\ba{\begin{eqnarray}}
\def\ea{\end{eqnarray}}

\def\o{\omega}
\def\Re{{\rm Re}}

\def\det{{\rm det}}

\def\log{\,{\rm log}\,}

\def\o{\omega}

\def\e{\varepsilon}

\def\o{\omega}

\def\na{\nabla}

\def\p{\prod}

\def\na{{\nabla}}

\def\[{{\bf [}}
\def\]{{\bf ]}}

\def\p{\partial}

\begin{document}
\starttext \baselineskip=15pt \setcounter{footnote}{0}
\newtheorem{theorem}{Theorem}
\newtheorem{lemma}{Lemma}
\newtheorem{definition}{Definition}
\newtheorem{proposition}{Proposition}

\begin{center}
{\Large \bf  THE ANOMALY FLOW ON UNIMODULAR LIE GROUPS
\footnote{Work supported in part by the National Science Foundation under NSF Grant DMS-12-66033.}}
\bigskip\bigskip

\centerline{Duong H. Phong, Sebastien Picard,
and Xiangwen Zhang}

\end{center}

\medskip

\begin{abstract}

{\small The Hull-Strominger system for supersymmetric vacua of the heterotic string allows general unitary Hermitian connections with torsion and not just the Chern unitary connection. Solutions on unimodular Lie groups exploiting this flexibility were found by T. Fei and S.T. Yau. The Anomaly flow is a flow whose stationary points are precisely the solutions of the Hull-Strominger system. Here we examine its long-time behavior on unimodular Lie groups with general unitary Hermitian connections. We find a diverse and intricate behavior, which depends very much on the Lie group and the initial data.}

\end{abstract}

\section{Introduction}
\setcounter{equation}{0}

The Hull-Strominger system \cite{H1,H2, S} is a system of equations for supersymmetric vacua of the heterotic string which generalizes the well-known Calabi-Yau compactifications found by Candelas, Horowitz, Strominger, and Witten \cite{CHSW}. It is also of considerable interest from the point of view of non-K\"ahler geometry and non-linear partial differential equations
\cite{FY1, FY2, PPZ1, PPZ3, PPZ4}. While more special solutions continue
to be constructed (see e.g. \cite{Fe1, Fe2,FHP, FTY, FeY, FIUV, OUV, FY1,FY2, GF, LY} and references therein), the complete solution of the Hull-Strominger system still seems distant. In \cite{PPZ2}, a general strategy was proposed, namely to look for solutions as stationary points of a flow, called there the Anomaly flow. Despite its unusual original formulation as a flow of $(2,2)$-forms, the Anomaly flow turns out to have some remarkable properties, including some suggestive analogies with the Ricci flow \cite{PPZ5}, and the fact \cite{PPZ6} that it can recover the famous solutions found in 2006 by Fu and Yau \cite{FY1, FY2}. Even so, very little is known at this time about the Anomaly flow in general, and it is important to gain insight from more examples.

\medskip
Part of the interest in the Hull-Strominger system lies in it allowing metrics on complex manifolds which have non-vanishing torsion. For such metrics, there is actually a whole line of natural unitary connections, namely the Yano-Gauduchon line \cite{G} passing by the Chern unitary connection and the Bismut connection. 
The Hull-Strominger system can be formulated with any specific choice along this line \cite{H2}, and not just the Chern unitary connection. This flexibility is essential in the solution found by Fei and Yau on unimodular Lie groups \cite{FeY}, using an ansatz originating from \cite{AGF} and \cite{BM}. The advantage of unimodular Lie groups is that any left-invariant Hermitian metric is automatically balanced, so if we fix a holomorphic vector bundle with a Hermitian-Yang-Mills invariant metric, the Hull-Strominger will reduce to a single equation. A difficulty in using general Hermitian connections is that their Riemannian curvature tensor $Rm$ will in general have components of all $(2,0)$, $(1,1)$, and $(0,2)$ types. But Fei and Yau discovered that, for unimodular Lie groups, the term ${\rm Tr}(Rm\wedge Rm)$ reduces to a $(2,2)$-form, and that solutions to the Hull-Strominger system can be found with any connection on the Yano-Gauduchon line, except for the Chern unitary connection and the Lichnerowicz connection.

\medskip
The purpose of this paper is analyze the Anomaly flow on unimodular Lie groups. As had been stressed earlier, the study of the Anomaly flow has barely begun, and there is a strong need for good examples. Even in the case of toric fibrations analyzed in \cite{PPZ6}, the convergence of the flow was established only for a particular type of initial data, and a fuller understanding of the long-time behavior of the flow is still not available. Unlike the more familiar cases of the K\"ahler-Ricci and the Donaldson heat flows, which are known to converge whenever a stationary point exists, the Anomaly flow is expected to behave differently depending on the initial data. This can be traced back to the fact that the conformally balanced condition is weaker than the K\"ahler condition, and that the positivity and closedness of a $(2,2)$-form seem to carry less information than the positivity and closedness of a $(1,1)$-form. The case of Lie groups is particular appealing, as the flow reduces then to a system of ordinary differential equations, and its formalism can be readily extended to general 
connections on the Yano-Gauduchon line. Our main results are as follows.

\medskip
Let $X$ be a complex $3$-dimensional Lie group. Fix a basis $\{e_a\}$ of left-invariant holomorphic vector fields on $X$, and let $\{e^a\}$ be the dual basis of holomorphic forms. 
Let $c_{ab}^d$ denote the structure constants of $X$ defined by this basis
\bea
\label{structureconstants}
[e_a,e_b]=\sum_d e_d\,c^d{}_{ab}.
\eea
The Lie group $X$ is said to be unimodular if the structure constants satisfy the condition
\bea
\label{unimodular}
\sum_a c^a{}_{ab}=0
\eea
for any $b$. This condition is invariant under a change of basis of holomorphic vector fields (see e.g. Section \S \ref{unimod-section} for a proof). The basis $\{e_a\}$ also defines a holomorphic, nowhere vanishing $(3,0)$-form on $X$, given by
\bea
\label{Omega}
\Omega=e^1\wedge e^2\wedge e^3.
\eea
We define the norm of $\Omega$ with respect to any metric $\o$ by
\bea
\label{NormOmega}
\|\Omega\|^2
=i\Omega\wedge \bar\Omega\, \bigg( {\o^3\over 3!} \bigg)^{-1}.
\eea

\medskip
\begin{theorem}
\label{Main1}
Let $X$ be a unimodular complex $3$-dimensional Lie group.
Let $t\to \o(t)$ be the Anomaly flow on $X$, as defined in (\ref{AnomalyFlow}) below. Set $\o(t)
=ig_{\bar ab}(t)\,e^b\wedge\bar e^a$. Then the Anomaly flow is given explicitly by
\bea \label{flowmetric}
\p_tg_{\bar{a} b} &=&
{1\over 2\|\Omega\|_{\o}} g^{d \bar{p}}  \bigg( g_{\bar{i} s}  \overline{c^i{}_{a p}}  c^s{}_{bd} - {\alpha' \tau \over 4} \, g_{\bar{\ell} i}  g^{n \bar{j}} \, \overline{c^{m}{}_{ap}} \overline{c^{\ell}{}_{mj}} c^i{}_{s n} c^s{}_{bd}\bigg).
\eea
Here $\alpha'$ is the string tension, and $\tau=2\kappa^2(2\kappa-1)$, where $\kappa$ indicates the connection (\ref{connection_rules}).
\end{theorem}

\medskip

In three dimensions, as listed in Fei and Yau \cite{FeY} and Knapp \cite{Kn}, the unimodular complex Lie groups consist precisely of those whose Lie algebra is isomorphic to that of either the Abelian Lie group ${\bf C}^3$, the nilpotent Heisenberg group, the solvable group of (complexfications of) rigid motions in ${\bf R}^2$, or the semi-simple group $SL(2,{\bf C})$. The Anomaly flow behaves differently in each case. Fixing the structure constants in each case as spelled out in \S 3, we have:
\begin{theorem}
\label{Main2} Assume that $\alpha'\tau>0$.

{\rm (a)} When $X={\bf C}^3$, any metric is a stationary point for the flow, and the flow is consequently stationary for any initial metric $\o(0)$.

{\rm (b)} When $X$ is nilpotent, there is no stationary point. Consequently the flow cannot converge for any initial metric. If the initial metric is diagonal, then the metric remains diagonal along the flow, the lowest eigenvalue is constant, while the other two eigenvalues tend to $+\infty$ at a constant rate.

{\rm (c)} When $X$ is solvable, the stationary points of the flow are precisely the metrics with 
\bea
g_{\bar 12}=\overline{g_{\bar 21}}=0,
\qquad {\alpha'\tau\over 4}g^{3\bar 3}=1.
\eea 
The Anomaly flow is asymptotically instable near any stationary point. However,
the condition $g_{\bar 12}=\overline{g_{\bar 12}}=0$ is preserved along the flow, and for any initial metric satisfying this condition, the flow converges to a stationary point.

{\rm (d)} When $X=SL(2,{\bf C})$, there is a unique stationary point, given by the diagonal metric
\bea
g_{\bar ab}={\alpha'\tau\over 2}\,\delta_{ab}.
\eea
The linearization of the flow at the fixed point admits both positive and negative eigenvalues. In particular, the flow is asymptotically instable. 

\end{theorem}

\medskip
The paper is organized as follows. In \S 2, the terms in the main equation in the Hull-Strominger system, namely the anomaly cancellation condition, are worked out in detail for connections on the Yano-Gauduchon line for left-invariant metrics. While the arguments are directly inspired by those of Fei and Yau, we have adopted a more explicit component formalism that should make the proof as well as the resulting formulas more accessible and convenient for geometric flows. The Anomaly flow is analyzed in \S 3. Theorem \ref{Main1} is proved in \S 3.2,
and Theorem \ref{Main2} is proved in \S 3.3, beginning in each case with the identification of the stationary points, and ending with the analysis of the Anomaly flow proper.

\section{The curvature of unitary connections
with torsion on Lie groups}
\setcounter{equation}{0}

In this section, we derive formulas for the curvature of an arbitrary connection on the Yano-Gauduchon line for Lie groups. They have been obtained before by Fei and Yau \cite{FeY}, but we need a formulation convenient for our subsequent study of the Anomaly flow.

\subsection{Unitary connections with torsion}

Let $(X, \o)$ be a complex manifold of complex dimension $n$, equipped with a Hermitian metric $\o$. We shall denote by $z^j$, $1\leq j\leq n$, local holomorphic coordinates on $X$, and by $\xi^\alpha$, $1\leq\alpha\leq 2n$, real smooth local coordinates on $X$. In particular $g_{\bar kj}$ denotes the metric in holomorphic coordinates, and $g_{\alpha\beta}$ the metric in real coordinates. We consider only unitary connections on $X$, that is connections $\na$ on $T^{1,0}(X)$ which preserve the metric, $\na g_{\bar kj}=0$. Our notation for covariant derivatives is
\bea
\label{covariantderivatives1}
\na_\alpha V^\gamma=\p_\alpha V^\gamma+A^\gamma{}_{\alpha\beta}V^\beta.
\eea
and our conventions for the curvature $R_{\alpha\beta}{}^\gamma{}_\delta$ and the torsion tensors $T^\gamma{}_{\alpha\beta}$ are,
\bea
\label{curvature0}
[\na_\beta,\na_\alpha]V^\gamma
=
R_{\alpha\beta}{}^\gamma{}_\delta V^\delta
+
T^\delta{}_{\alpha\beta}\na_\delta V^\gamma,
\eea
with similar formulas when using the complex coordinates $z^j, \bar z^j$, $1\leq j\leq n$. For example, in complex coordinates, the components of the torsion are
\bea
\label{torsion}
&&
T^p{}_{\bar jm}=A_{\bar jm}^p,
\qquad
T^{\bar p}{}_{\bar jm}=-A^{\bar p}_{m\bar j},
\nonumber\\
&&
T^p{}_{jm}=A_{jm}^p-A_{mj}^p,
\qquad
T^{\bar p}{}_{\bar j\bar m}
=
A_{\bar j\bar m}^{\bar p}-A_{\bar m\bar j}^{\bar p}.
\eea

\smallskip

Among the unitary connections, of particular interest are the following three special connections:

\medskip
(1) The Levi-Civita connection $\na^L$ on $T(X)$, characterized by the fact that its torsion is $0$. In general, it does not preserve the complex structure, and in particular, it does not induce a connection on the complex tangent space $T^{1,0}(X)$. It preserves the complex structure $J$ if and only if $\o$ is K\"ahler.

\smallskip
(2) The Chern connection $\na^C$, characterized by the fact that it preserves the complex structure, and that $\na_{\bar j}^CW^k=\p_{\bar j}W^k$, i.e.,
$A_{\bar j m}^{p}=0$.  By unitarity, it follows that
$A_{jm}^p=g^{p\bar q}\p_j g_{\bar qm}$. Thus the components of the torsion of the Chern connection in the first line of (\ref{torsion}) vanish, and the torsion reduces to a section of the following bundle
\bea
{\cal T}=T+\bar T\in (\Lambda^{2,0}\otimes T^{1,0})\oplus (\Lambda^{0,2}\otimes T^{0,1}). 
\eea
with 
\bea
&&
T={1\over 2}\,T^\ell{}_{jk}\,{\p\over\p z^\ell}\otimes (dz^k\wedge dz^j),
\quad
\bar T
={1\over 2}\,\bar T^{\bar\ell}{}_{\bar j\bar k}
\,{\p\over\p \bar z^\ell}\otimes (d\bar z^k\wedge d\bar z^j)
\quad
\bar T^{\bar\ell}{}_{\bar j\bar k}
=
\overline{T^\ell{}_{jk}}.\nonumber
\eea
The torsion of the Chern connection is $0$ exactly when $\o$ is K\"ahler.

\smallskip
(3) The Bismut connection $\na^{\cal H}$. Connections obtained by modifying the Levi-Civita connection by a $3$-form seem to have been written down first by Yano in his book \cite{Ya}, which became the standard reference for it in the physics literature on supersymmetric sigma models (see e.g. \cite{GHR, HW, H1, S} and references therein). They were subsequently rediscovered by Bismut \cite{B}, 
who also identified among them the connection preserving the complex structure,
which is now known as the Bismut connection.
Thus the Bismut connection can be written 
in two different ways, depending on whether we use $\na^L$ or $\na^C$ as reference connection.
To write it with $\na^L$ as reference connection, introduce the following real $3$-form
\footnote{It is important to distinguish the form $H\in\Lambda^{2,1}$ from the torsion tensor $T\in\Lambda^{2,0}\otimes T^{1,0}$, although their coefficients are the same, and they can be recovered from one another. Perhaps a good analogue is the metric $g_{\bar kj}$ and the corresponding symplectic form $\o=ig_{\bar kj}\,dz^j\wedge d\bar z^k$, which are distinct objects, although they can be recovered from one another, once the complex structure is fixed.},
\bea
{\cal H}
={1\over 4}(H+\bar H)
\in \Lambda^{2,1}\oplus\Lambda^{1,2}\subset\Lambda^3
\eea
where
\bea
&&
H={1\over 2}\,g_{\bar m\ell}\,T^\ell{}_{jk}\,dz^k\wedge dz^j\wedge d\bar z^m,
\qquad
\bar H
={1\over 2}\,g_{\bar \ell m}\,\bar T^{\bar \ell}{}_{\bar j\bar k}\,d\bar z^k\wedge d\bar z^j\wedge dz^m,
\nonumber\\
&&
H_{\bar m jk}=g_{\bar m\ell}T^\ell{}_{jk},
\qquad
\bar H_{m\bar j\bar k}=g_{\bar\ell m}\bar T^\ell{}_{\bar j\bar k}.
\eea
Note that, in terms of the Hermitian form $\o=ig_{\bar kj}dz^j\wedge d\bar z^k$, we have
\bea
\label{H-omega}
H=i\p\o,\qquad
\bar H=-i\bar\p\o.
\eea
Next, introduce the $1$-form $S^{\cal H}=(S_\alpha^{\cal H})d\xi^\alpha$ valued in the space of endomorphisms anti-symmetric with respect to the metric $g_{\alpha\beta}$ by
setting
\footnote{An endomorphism $S^\beta{}_\gamma$ is antisymmetric with respect to a metric $g_{\rho\beta}$ if $\<SW,U\>=-\<W,SU\>$, which means that $g_{\beta\rho}S^\beta{}_\gamma=-S^\mu{}_{\rho}g_{\gamma\mu}$, which means that $S_{\rho\gamma}$ is anti-symmetric, if we raise and power indices using the metric $g_{\beta\rho}$.}
\bea
&&
{\cal H}\equiv {1\over 3!}\,{\cal H}_{\alpha\beta\gamma}\,d\xi^\gamma\wedge d\xi^\beta\wedge d\xi^\alpha,
\nonumber\\
&&
(S_\alpha^{\cal H})^\beta{}_\gamma
=
2{\cal H}_{\alpha\rho\gamma}\,g^{\rho\beta}.
\eea
The Bismut connection is then defined by
\bea
\label{YB1}
\na^{\cal H}=\na^L+S^{\cal H}.
\eea
This formulation shows that it is unitary and that its torsion can be identified with a $3$-form.
If we use the Chern connection as the reference connection, we can also write
\bea
\label{YB2}
\na_j^{\cal H}W^p=\na_j^C W^p-T^p{}_{jk}W^k,
\qquad
\na_{\bar j}^{\cal H}W^p
=
\na_{\bar j}^CW^p+g^{p\bar q}g_{\bar mr}\bar T^{\bar m}{}_{\bar j\bar q}W^r.
\eea
This formulation shows that the Bismut connection
is unitary and that it manifestly preserves the complex structure.
The equivalence of (\ref{YB1}) and (\ref{YB2}) can be found in \cite{B}.

\bigskip

The two connections $\na^C$ and $\na^{\cal H}$ determine a line of connections $\nabla^{(\kappa)}$ which are all unitary and preserve the complex structure. Indeed, in terms of an orthonormal frame, a unitary connection is characterized by a Hermitian matrix, and a real linear combination of Hermitian matrices is again Hermitian. Also a linear combination of connections on $T^{1,0}(X)$ is again a connection on $T^{1,0}(X)$, hence our assertion. We refer to this line of connections as the Yano-Gauduchon line of canonical connections \cite{G}. With $\na^C$ as reference connection, the connection $\na^{(\kappa)}$ can be expressed as
\bea \label{connection_rules}
\na_j^{(\kappa)}W^p=\na_j^C W^p-\kappa\,T^p{}_{jk}W^k,
\qquad
\na_{\bar j}^{(\kappa)}W^p
=
\na_{\bar j}^CW^p+\kappa\,g^{p\bar q}g_{\bar mr}\bar T^{\bar m}{}_{\bar j\bar q}W^r.
\eea
Clearly $\kappa=0$ corresponds to the Chern connection and $\kappa=1$ to the Bismut connection. The case $\kappa=1/2$ is known as the Lichnerowicz connection.

\subsection{Evaluation of curvature forms}

Let $(X,\o)$ be now a complex Lie group with left-invariant metric $\o$. Let $e_1,\dots , e_n \in {\bf g}$ be an orthonormal frame of left-invariant holomorphic vector fields on $X$, and let $e^1,\dots , e^n$ be the dual frame of holomorphic 1-forms, so that
\be
\o = i \sum_a e^a \wedge \bar{e}^a.
\ee
The structure constants of the Lie algebra ${\bf g}$ in this basis are defined in
(\ref{structureconstants}) and denoted
by $c^d{}_{ab}$. They satisfy the Jacobi identity
\be \label{jacobi-id}
c^q{}_{ir}c^r{}_{jk} + c^q{}_{kr}c^r{}_{ij} + c^q{}_{jr}c^r{}_{ki} =0.
\ee
By Cartan's formula for the exterior derivative, we have
\be \label{Maurer-Cartan}
\p e^a = {1 \over 2} c^a{}_{bd} \, e^d \wedge e^b
\ee
which is the Maurer-Cartan equation. 

\medskip
Our goal is to determine the torsion and curvature of an arbitrary connection $\na^{(\kappa)}$ on the Yano-Gauduchon line in terms of the frame $e^1,\cdots,e^n$ and the structure constants $c^d{}_{ab}$. First, we express the connection forms for $\na^{(\kappa)}$ in terms of the structure constants of $X$:

\begin{lemma}
Let $\na^{(\kappa)}$ be a connection on the Yano-Gauduchon line for a left-invariant Hermitian metric $\o$ on the Lie group $X$. Let $e^1,\cdots,e^n$ be a left invariant orthonormal basis for $\o$, and let $c^d{}_{ab}$ be the corresponding structure constants, as defined above. Then for any vector field $V=V^a e_a$,
\be \label{Lie-groups-connection}
\na_b^{(\kappa)} V^a = \p_{e_b} V^a + \kappa c^a{}_{bd} V^d, \ \ \na_{\bar{b}}^{(\kappa)} V^a = \p_{\bar{e}_b} V^a - \kappa \overline{c^d{}_{ba}} V^d.
\ee
\end{lemma}

\medskip
\noindent
{\it Proof.} We use the expression for $\na^{(\kappa)}$ with the Chern connection $\na^C$ as reference connection. For this, we need the torsion form $H$ of the Chern connection.
Taking the exterior derivative of $\o$ and applying (\ref{Maurer-Cartan}) gives
\be
\label{Homega}
H=i \p \o = - {1 \over 2} \sum c^d{}_{ab} \, e^b \wedge e^a \wedge \bar{e}^d.
\ee
Therefore
\be
H_{\bar d ab}=( i \p \o)_{\bar{d} ab} = - c^d{}_{ab}.
\ee
Next, since the coefficients of the metric $\o$ are constant in an orthonormal frame, the Chern connection reduces to the exterior derivative: $\na^C_j = \p_j$ in this frame. The lemma follows now from the formula (\ref{connection_rules}) giving $\na^{(\kappa)}$ in terms of $\na^C$, $H$, and the parameter $\kappa$.

\bigskip

Henceforth, we shall denote $\nabla^{(\kappa)}$ by just $\na$ for simplicity. We obtain
\be
A^a{}_{bd}= \kappa c^a{}_{bd}, \ \ A^a{}_{\bar{b} d}= - \kappa \overline{c^d{}_{ba}}, \ \ A^{\bar{a}}{}_{b \bar{d}}= - \kappa c^d{}_{ba}, \ \ A^{\bar{a}}{}_{\bar{b} \bar{d}}= \kappa \overline{c^a{}_{bd}}.
\ee

\medskip
\begin{lemma}
The components of the curvature of the connection $\na^{(\kappa)}$ in the orthonormal frame $e^1,\cdots,e^n$ are given by
\bea
\label{curvature}
&&
R_{kj}{}^p{}_q =(\kappa-\kappa^2) c^r{}_{kj} c^p{}_{rq},
\qquad
R_{\bar{k} \bar{j}}{}^p{}_q = - (\kappa-\kappa^2) \overline{ c^r{}_{kj} c^q{}_{rp}}
\nonumber\\
&&
R_{\bar{k} j}{}^p{}_q = \kappa^2 (-c^p{}_{jr} \overline{c^q{}_{kr}}  + \overline{c^r{}_{kp}} c^r{}_{jq}).
\eea
\end{lemma}

\medskip
\noindent
{\it Proof.} The defining formula (\ref{curvature0}) for the curvature
and the torsion $t^d{}_{ba}$ of the connection $\na^{(\kappa)}$ in coordinates
becomes the following formula
in the frame $\{ e_a \}$,
\be \label{curv-frame-2}
[\na_a, \na_b] W^c = R_{ba}{}^c{}_d W^d + t^d{}_{ba} \na_d W^c + c^d{}_{ab} \na_d W^c
\ee
where the last term on the right hand side is the contribution of the commutator $[e_a,e_b]$. H
We now compute the left hand side
using the formulas of the previous lemma. We find
\bea
 R_{kj}{}^p{}_q W^q + t^q{}_{kj} \na_q W^p + c^q{}_{jk} \na_q W^p
&=&
[\na_j, \na_k] W^p\nonumber\\
&=& c^q{}_{jk} \p_{e_q} W^p+
\kappa^2 (-c^p{}_{kr} c^r{}_{js}+c^p{}_{jr}c^r{}_{ks})W^s
\nonumber\\
&&
+2\kappa c^r{}_{kj} \na_rW^p.
\eea
By (\ref{Lie-groups-connection}) and the Jacobi identity (\ref{jacobi-id}), we see that the $(2,0)$ component is
\bea\label{2,0R}
R_{kj}{}^p{}_q =- \kappa c^r{}_{jk} c^p{}_{rq} + \kappa^2 (-c^p{}_{kr} c^r{}_{jq}+c^p{}_{jr}c^r{}_{kq}) = (\kappa-\kappa^2) c^r{}_{kj} c^p{}_{rq}.
\eea
This proves the first formula in the lemma.
Next, we have
\bea
R_{\bar{k} \bar{j}}{}^p{}_q W^q + \bar{t}^{\bar{q}}{}_{\bar{k} \bar{j}} \na_{\bar{q}} W^p + \overline{c^q{}_{jk}} \na_{\bar{q}} W^p
&=&
[\na_{\bar{j}}, \na_{\bar{k}}] W^p \nonumber\\
&=&
\overline{c^q{}_{jk}} \p_{\bar{e}_q} W^p + \kappa^2 (\overline{ c^r{}_{jp} c^s{}_{kr} - c^r{}_{kp} c^s{}_{jr}})W^s
\nonumber\\
&&
- 2 \kappa \overline{c^r{}_{kj}} \na_{\bar r}W^p.
\eea
We see that the $(0,2)$ component is
\be
R_{\bar{k} \bar{j}}{}^p{}_q = \kappa \overline{c^r{}_{jk}} \overline{c^q{}_{rp}} + \kappa^2 (\overline{ c^r{}_{jp} c^q{}_{kr} - c^r{}_{kp} c^q{}_{jr}}) = - (\kappa-\kappa^2) \overline{ c^r{}_{kj} c^q{}_{rp}}.
\ee
which is the second formula in the lemma. Finally,
we have
\bea 
R_{\bar{k} j}{}^p{}_q W^q - A^{\bar{q}}{}_{j \bar{k}} \na_{\bar{q}} W^p + A^{q}{}_{\bar{k} j} \na_{q} W^p
&=& [\na_j, \na_{\bar{k}}] W^p \nonumber\\
&=&
 \kappa c^k{}_{jr} \na_{\bar{r}} W^p
- \kappa \overline{c^j{}_{kr}} \na_r W^p
\nonumber\\
&&
+\kappa^2 (-c^p{}_{jr} \overline{c^q{}_{kr}}  + \overline{c^r{}_{kp}} c^r{}_{jq})W^q.
\eea
We see that the $(1,1)$ component is
\be
R_{\bar{k} j}{}^p{}_q = \kappa^2 (-c^p{}_{jr} \overline{c^q{}_{kr}}  + \overline{c^r{}_{kp}} c^r{}_{jq}).
\ee
This completes the proof of the lemma.

\bigskip

We note our expressions satisfy the usual symmetries
\bea
\label{Riem-kj-sym}
&&
R_{kj}{}^p{}_s = - R_{jk}{}^p{}_s, \ \ R_{\bar k\bar j}{}^p{}_s = - R_{\bar j\bar k}{}^p{}_s
\\
&&
\label{Riem-kj-sym1}
R_{\bar k \bar j}{}^p {}_q = - \overline{R_{kj}{}^q{}_p}, \ \ R_{\bar k j}{}^p{}_q = \overline{R_{\bar j k}{}^q{}_p}.
\eea

\medskip

Next, we compute ${\rm Tr}\, (Rm \wedge Rm)$ on the Lie group $X$. This computation was first done by Fei-Yau \cite{FeY}. Their computation showed in particular that ${\rm Tr}\, (Rm \wedge Rm)$ is a $(2,2)$ form along the Yano-Gauduchon line. Here we recover their result using our formalism.

\begin{lemma}
Let $\na^{(\kappa)}$ be an arbitrary connection on the Yano-Gauduchon line, as before. Assume that ${\rm dim}\,X=3$. Then ${\rm Tr}\, (Rm \wedge Rm)$ is a $(2,2)$-form, given explicitly by
\bea \label{1,1-comp}
{\rm Tr}\, (Rm \wedge Rm)_{\bar{k} \bar{\ell} i j}
&=& \tau \, \overline{c^r{}_{k \ell} c^s{}_{rp}} c^q{}_{ij} c^s{}_{qp}.
\eea
where we have set as before
$\tau= 2 \kappa^2 (2\kappa-1)$. 
\end{lemma}

\medskip
\noindent
{\it Proof.} Because $X$ is assumed to have dimension $3$, the $(4,0)$ and $(0,4)$ components of ${\rm Tr}(Rm\wedge Rm)$ are automatically $0$.
Next, we show the vanishing of the $(3,1)$-component,
\be \label{3,1-Lie0}
{\rm Tr}\, (Rm \wedge Rm)_{\bar{\ell} ijk} = 2( R_{ij}{}^p{}_s R_{\bar{\ell} k}{}^s{}_p + R_{ki}{}^p{}_s R_{\bar{\ell} j}{}^s{}_p  + R_{jk}{}^p{}_s R_{\bar{\ell} i}{}^s{}_p).
\ee
Fixing an index $(ijk,\ell)$, we use the previously computed formulas for the curvature of a Lie group to obtain
\be
R_{ij}{}^p{}_s R_{\bar{\ell} k}{}^s{}_p =  (\kappa-\kappa^2) \kappa^2c^r{}_{ij} c^p{}_{rs}(-  c^s{}_{kq} \overline{c^p{}_{\ell q}} +  \overline{c^q{}_{\ell s}} c^q{}_{kp}).
\ee
Applying the Jacobi identity (\ref{jacobi-id}),
\bea
& \ & R_{ij}{}^p{}_s R_{\bar{\ell} k}{}^s{}_p \\
&=&  \kappa^3(1-\kappa) ( -c^p{}_{ir} c^r{}_{js}c^s{}_{kq} \overline{c^p{}_{\ell q}} +  c^p{}_{jr} c^r{}_{is} c^s{}_{kq} \overline{c^p{}_{\ell q}} - c^q{}_{kp} c^p{}_{jr} c^r{}_{is} \overline{c^q{}_{\ell s}} +  c^q{}_{kp} c^p{}_{ir} c^r{}_{js} \overline{c^q{}_{\ell s}})  \nonumber\\
&=& 
 \kappa^3(1-\kappa) ( - c^p{}_{ir} c^r{}_{js} c^s{}_{kq} + c^p{}_{jr} c^r{}_{is} c^s{}_{kq} - c^p{}_{ks} c^s{}_{jr} c^r{}_{iq} + c^p{}_{ks} c^s{}_{ir} c^r{}_{jq} ) \overline{ c^p{}_{\ell q}}.
\eea
If we denote $F_{ijk}{}^p{}_{q} = c^p{}_{ir} c^r{}_{js} c^s{}_{kq}$, then
\bea
R_{ij}{}^p{}_s R_{\bar{\ell} k}{}^s{}_p =   \kappa^3(1-\kappa) ( - F_{ijk}{}^p{}_{q} + F_{jik}{}^p{}_q - F_{kji}{}^p{}_q + F_{kij}{}^p{}_q)\, \overline{c^p{}_{\ell q}}.
\eea
Upon cyclically permuting $(ijk)$, we see that
\be
R_{ij}{}^p{}_s R_{\bar{\ell} k}{}^s{}_p + R_{ki}{}^p{}_s R_{\bar{\ell} j}{}^s{}_p  + R_{jk}{}^p{}_s R_{\bar{\ell} i}{}^s{}_p =0.
\ee
Therefore
\be
{\rm Tr}\, (Rm \wedge Rm)_{\bar{\ell} ijk} = 0
\ee
as claimed.
Next, we compute the $(1,3)$ component,
\bea \label{1,3-Lie1}
{\rm Tr}\, (Rm \wedge Rm)_{\bar{i} \bar{j} \bar{k} \ell} &=& 2( R_{\bar{i} \bar{j}}{}^p{}_s R_{\bar{k} \ell}{}^s{}_p + R_{\bar{j} \bar{k}}{}^p{}_s R_{\bar{i} \ell}{}^s{}_p +  R_{\bar{k} \bar{i}}{}^p{}_s R_{\bar{j} \ell}{}^s{}_p).
\eea
Using the symmetric property (\ref{Riem-kj-sym1}) and the vanishing of the $(3,1)$ part, we see that the $(1,3)$ part also vanishes.

\medskip
Finally, we evaluate the $(2,2)$-component,
\be \label{2,2-Lie1}
{\rm Tr}\, (Rm \wedge Rm)_{\bar{k} \bar{\ell} i j} = 2 R_{\bar{k} \bar{\ell}}{}^p{}_s R_{ij}{}^s{}_p + 2(R_{\bar{k} j}{}^p{}_s R_{\bar{\ell} i}{}^s{}_p - R_{\bar{k} i}{}^p{}_s R_{\bar{\ell} j}{}^s{}_p) .
\ee
From the previously established formulas, we have
\bea
R_{\bar{k} \bar{\ell}}{}^p{}_s R_{ij}{}^s{}_p = - (\kappa-\kappa^2)^2\,  \overline{c^r{}_{k\ell} c^s{}_{rp}} c^q{}_{ij} c^s{}_{qp}.
\eea
Next,
\bea \label{1,1-1,1}
 R_{\bar{k} j}{}^p{}_s R_{\bar{\ell} i}{}^s{}_p&=& \kappa^4(- c^p{}_{jr} \overline{c^s{}_{kr}} + \overline{c^r{}_{kp}} c^r{}_{js}) ( -c^s{}_{iq} \overline{c^p{}_{\ell q}} + \overline{c^q{}_{\ell s}} c^q{}_{ip}) \\
&=& \kappa^4 ( c^p{}_{jr} \overline{c^s{}_{kr}}c^s{}_{iq} \overline{c^p{}_{\ell q}} + c^q{}_{ip} \overline{c^r{}_{kp}} c^r{}_{js}\overline{c^q{}_{\ell s}}  -  c^p{}_{jr} \overline{c^s{}_{kr}} \overline{c^q{}_{\ell s}} c^q{}_{ip}  - c^s{}_{iq}\overline{c^p{}_{\ell q}}  \overline{c^r{}_{kp}}  c^r{}_{js}). \nonumber
\eea
Some cancellation occurs and we are left with
\bea
& \ & (R_{\bar{k} j}{}^p{}_s R_{\bar{\ell} i}{}^s{}_p - R_{\bar{k} i}{}^p{}_s R_{\bar{\ell} j}{}^s{}_p) \nonumber\\
&=&  \kappa^4 (- \overline{c^s{}_{kr}} \overline{c^q{}_{\ell s}} c^q{}_{ip} c^p{}_{jr}   - \overline{c^p{}_{\ell q}}  \overline{c^r{}_{kp}}  c^r{}_{js}c^s{}_{iq} +   \overline{c^s{}_{kr}} \overline{c^q{}_{\ell s}} c^q{}_{jp} c^p{}_{ir} +\overline{c^p{}_{\ell q}}\overline{c^r{}_{kp}}  c^r{}_{is} c^s{}_{jq})\nonumber\\
&=& \kappa^4 \left( \overline{c^s{}_{kr}} \overline{c^q{}_{\ell s}} (-c^q{}_{ip} c^p{}_{jr} + c^q{}_{jp} c^p{}_{ir}) + \overline{c^p{}_{\ell q}}\overline{c^r{}_{kp}} (-c^r{}_{js} c^s{}_{iq} + c^r{}_{is} c^s{}_{jq})\right)\nonumber
\\
&=& \kappa^4 (\overline{c^s{}_{kr}} \overline{c^q{}_{\ell s}} c^q{}_{rp} c^p{}_{ij} -\overline{c^p{}_{\ell q}}\overline{c^r{}_{kp}} c^r{}_{qs} c^s{}_{ij})\nonumber
\\
&=& \kappa^4 (\overline{c^s{}_{kr}} \overline{c^q{}_{\ell s}} - \overline{c^s{}_{\ell r} c^q{}_{ks}}  ) c^q{}_{rp} c^p{}_{ij}\nonumber
\\
&=& \kappa^4 \overline{c^q{}_{rs} c^s{}_{k\ell}} c^q{}_{rp} c^p{}_{ij}.
\eea
Adding these two equations together, the terms of order $\kappa^4$ cancel and we obtain the desired formula.

\bigskip
Recall that we view the Lie algebra of $X$ as generated by a given basis of left-invariant holomorphic vector fields $e_a$, with structure constants $c^d{}_{ab}$. So far we have considered only the metric $\o=i\sum e^a\wedge \bar e^a$ defined by the condition that $e_a$ be orthonormal. We consider now the general left-invariant metric given by
\be 
\label{generalmetric}
\o = \sum\, g_{\bar{b} a} \, i e^a \wedge \bar{e}^b,
\ee
where $g_{\bar{b} a}$ is a positive-definite Hermitian matrix. We have 

\medskip
\begin{lemma}
\label{Lemma4}
Let $X$ be a 3-dimensional complex Lie group, with a given basis of left-invariant holomorphic vector fields $e_a$ with structure constants $c^d{}_{ab}$. Let $\o$ be the metric given by (\ref{generalmetric}), and let $\na^{(\kappa)}$ be the Hermitian connection on the Gauduchon line with parameter $\kappa$. Then
\be \label{TrR-wedge-R}
{\rm Tr}\, (Rm \wedge Rm) = {\tau \over 4} \, g_{\bar{\ell} i}  g^{n \bar{j}} \, \overline{c^{m}{}_{ab}} \overline{c^{\ell}{}_{mj}} c^i{}_{s n} c^s{}_{cd} \, e^d \wedge e^c \wedge \bar{e}^{b} \wedge \bar{e}^{a}.
\ee
and
\be
  i \p \bar{\p} \o - {\alpha' \over 4} {\rm Tr}\, (Rm \wedge Rm) = {1 \over 4} \bigg( g_{\bar{i} s}  \overline{c^i{}_{a b}}  c^s{}_{cd} - {\alpha' \tau \over 4} \, g_{\bar{\ell} i}  g^{n \bar{j}} \, \overline{c^{m}{}_{ab}} \overline{c^{\ell}{}_{mj}} c^i{}_{s n} c^s{}_{cd}\bigg) \,  e^d \wedge  e^c \wedge \bar{e}^b \wedge \bar{e}^a.
\ee
\end{lemma}

\medskip
\noindent
{\it Proof.} These formulas for general metrics follow from the ones obtained earlier for the metric $i\sum_a e^a\wedge\bar e^a$ after performing a change of basis.
More specifically, we let $P$ be a matrix such that
\be
\bar{P}^{\bar{a}}{}_{\bar{p}} g_{\bar{a} b} P^b{}_{q} = \delta_{\bar{p} q}.
\ee
Therefore, denoting $g^{a \bar{b}}$ to be the inverse of $g_{\bar{b} a}$, we have
\be \label{metric-vs-P}
g_{\bar{b} a} = (\bar{P}^{-1})^{\bar{r}}{}_{\bar{b}} (P^{-1})^r{}_{a}, \ \ g^{a \bar{b}} = P^a{}_{r} \bar{P}^{\bar{b}}{}_{\bar{r}}.
\ee
We now perform a change of basis and define 
\be
f_i = e_r P^r{}_{i}, \ \ [f_i,f_j] = f_r k^r{}_{ij}.
\ee
The induced transformation laws are
\be
f_i = e_r P^r{}_{i}, \ \ f^i = (P^{-1})^i{}_r e^{r}, \ \ e_i = f_r (P^{-1})^r{}_{i}, \ \ e^i = P^i{}_r f^{r},
\ee
\be
\label{transformation}
k^\ell{}_{ij} = (P^{-1})^\ell{}_s P^r{}_i P^q{}_{j} c^s{}_{rq}.
\ee
By construction, $\o$ is diagonal in the basis $\{ f_a \}$. From (\ref{1,1-comp}), we can compute the curvature in the basis $\{ f^a \}$,
\bea
{\rm Tr}\, (Rm \wedge Rm) = {\tau \over 4}  \, \overline{k^r{}_{a b} k^s{}_{rp}} k^q{}_{cd} k^s{}_{qp} \, f^d \wedge f^c \wedge \bar{f}^{b} \wedge \bar{f}^{a}.
\eea
Using the above transformation laws, ${\rm Tr}\, (Rm \wedge Rm)$ can be rewritten in terms of the $\{ e^a \}$ basis. A straightforward computation gives the first formula in the lemma. To obtain the second formula, we begin by computing $i\p\bar\p\o$ in the model case where the metric is $\sum_a ie^a\wedge\bar e^a$. We had already found $i\p\o$ in (\ref{Homega}). Differentiating again gives
\be \label{iddb-omega}
i \p \bar{\p} (i\sum_a e^a\wedge\bar e^a) = {1 \over 4} \sum \overline{c^\ell{}_{ab}} \, c^\ell{}_{cd} \, e^d \wedge e^c \wedge \bar{e}^b \wedge \bar{e}^a.
\ee
Reverting to the general metric $\o$ given by (\ref{generalmetric}) and performing the same change of bases as before, we find
\be \label{iddb-omega-general}
i \p \bar{\p} \o = {1 \over 4} g_{\bar{i} s}  \overline{c^i{}_{a b}}  c^s{}_{cd} \,  e^d \wedge  e^c \wedge \bar{e}^b \wedge \bar{e}^a.
\ee
Combining this formula with the one found previously for ${\rm Tr}(Rm \wedge Rm)$, we obtain the second formula stated in the lemma. Q.E.D.

\section{Hull-Strominger systems and Anomaly flows}
\setcounter{equation}{0}

We come now to the study of the Anomaly flow on the complex Lie group $X$.

\subsection{The Hull-Strominger system on unimodular Lie groups} \label{unimod-section}

First we recall the Hull-Strominger system \cite{H1,H2,S}. Let $X$ be a $3$-dimensional complex manifold with a nowhere vanishing holomorphic $(3,0)$-form $\Omega$. The Hull-Strominger system is a system of equations for a Hermitian metric $\o$ on $X$ and a holomorphic vector bundle $E\to X$ equipped with a Hermitian metric $H_{\bar\alpha\beta}$ satisfying
\bea
\label{HS}
&&
F^{0,2}=F^{0,2}=0,\qquad \o^2\wedge F^{1,1}=0
\nonumber\\
&&
i\p\bar\p\o-{\alpha'\over 4}
{\rm Tr}(Rm \wedge Rm -F\wedge F)=0
\nonumber\\
&&
d^\dagger\o=i(\p-\bar\p)\log \|\Omega\|_\o
\eea
where $F^{p,q}$ are the components of the Chern curvature $F\in \Lambda^2
\otimes End(E)$ of the metric $H_{\bar\alpha\beta}$, and $\|\Omega\|_\o$ denotes the norm of $\Omega$ with respect to the metric $\o$ as defined in (\ref{NormOmega}). The first equation in (\ref{HS}) is the familiar Hermitian-Yang-Mills equation, and its solution is well-known for given metric $\o$ by the theorem of Donaldson-Uhlenbeck-Yau \cite{D,UY}. Thus the most novel aspects in the Hull-Strominger system resides in the other two equations. It has been pointed out by Li and Yau \cite{LY} that the last equation is equivalent to the following condition of ``conformally balanced metric",
\bea
\label{balanced}
d(\|\Omega\|_\o\o^2)=0
\eea
which is a generalization of the balanced condition introduced in 1981 by Michelsohn \cite{M}.

\medskip
We shall take the bundle $E\to X$ to be trivial, $F$ to be $0$, and restrict ourselves to left-invariant metrics. In this case, the first equation in (\ref{HS}) is trivially satisfied. Furthermore, the norm $\|\Omega\|_\o$ is a constant function on $X$, and the conformally balanced condition reduces to exactly the balanced condition of Michelsohn \cite{M}, i.e., $d\omega^2=0$. 

\medskip
A key observation, also exploited earlier in the work of Fei and Yau \cite{FeY}, is that, on a unimodular complex Lie group, any left-invariant metric is unimodular.
This statement is well-known and to our knowledge first appeared in \cite{AG}. For the reader's convenience, we provide the brief argument. Recall that a complex Lie group is said to be unimodular if there exists a left-invariant basis of holomorphic vector fields with structure constants $c^d{}_{ab}$ satisfying
the condition (\ref{unimodular}).
In view of the transformation rule (\ref{transformation}) for structure constants under a change of basis of left-invariant vector fields, this statement holds for all bases if and only if it holds for some basis. Let now $\o$ be any invariant Hermitian metric on $X$, and express it in terms of any basis of holomorphic $e^a$ forms as (\ref{generalmetric}). A direct calculation using (\ref{Maurer-Cartan}) gives
\be
\p \o^2 = \sum {1 \over 2} (g_{\bar{p} r} g_{\bar{q} c} - g_{\bar{p} c} g_{\bar{q} r}) c^r{}_{ab} \, e^a \wedge e^b \wedge e^c \wedge \bar{e}^p \wedge \bar{e}^q.
\ee
Choosing $e^a$ to be orthonormal with respect to $\o$, we can assume that 
$g_{\bar{b} a} = \delta_{ba}$, and we 
readily see that the condition $d \o^2 = 0$ is equivalent to the unimodular condition
(\ref{unimodular}). 

\smallskip

Thus, on unimodular Lie groups, the Hull-Strominger system for a left-invariant metric reduces to the middle equation in (\ref{HS}).

\subsection{Proof of Theorem \ref{Main1}}

The Anomaly flow, introduced in \cite{PPZ2}, is a parabolic flow whose stationary points are solutions of the Hull-Strominger system. In the present setting, $X$ is a $3$-dimensional unimodular complex Lie group with a basis $\{e_a\}$ of left-invariant holomorphic vector fields, $\Omega=e^1\wedge e^2\wedge e^3$, and the bundle $E\to X$ is taken to be trivial. 
Then the Anomaly flow \cite{PPZ2} is defined to be the following flow of $(2,2)$-forms,
\bea
\label{AnomalyFlow}
\p_t(\|\Omega\|_\o\o^2)
=
i \p \bar{\p} \o - {\alpha' \over 4} {\rm Tr}\, (Rm \wedge Rm).
\eea
with any given initial data of the form $\|\Omega\|_{\o_0}\o_0^2$.

Theorem 1 of \cite{PPZ5} shows how to rewrite this flow as a curvature flow for the Hermitian metric $\o$. Setting 
$\Phi = i \p \bar{\p} \o - {\alpha' \over 4} {\rm Tr}\, (Rm \wedge Rm)$, we obtain
\be
\p_tg_{\bar{a} c} = -{1\over 2\|\Omega\|_{\o}} g^{d \bar{b}}\Phi_{\bar{a}\bar{b}cd}.
\ee
Substituting in the formulas obtained in Lemma \ref{Lemma4} for $\Phi$,
we obtain Theorem \ref{Main1}.

\subsection{Proof of Theorem \ref{Main2}}

We discuss the unimodular Lie group case by case, as listed in Theorem \ref{Main2}.

\subsubsection{The Abelian Case}

In this case, we have for all $a,b$,
\bea
[e_a,e_b]=0
\eea
and all the structure constants $c^d{}_{ab}$ vanish. The flow (\ref{flowmetric}) is static for all initial data, and part (a) of Theorem \ref{Main2} is immediate.

\subsubsection{Nilpotent case}
In this case, we may assume the Lie algebra satisfies the commutation relations 
\bea
[e_1,e_3]=0, \ [e_2,e_3]=0, \ [e_1,e_2] = e_3.
\eea 
It follows that $c^3{}_{12}=1$ and all other structure constants vanish. Substituting these structure constants into the flow (\ref{flowmetric}) gives the following system
\be\label{Nilpotentsystem}
\p_t g_{\bar{1} 1} = {g^{2 \bar{2}} g_{\bar{3} 3} \over 2 \| \Omega \|} , \ \ \p_t g_{\bar{1} 2} = - {g^{1 \bar{2}} g_{\bar{3} 3} \over 2 \| \Omega \|}, \ \ \p_t g_{\bar{2} 2} = {g^{1 \bar{1}} g_{\bar{3} 3} \over 2 \| \Omega \|}, \ \ \p_t g_{\bar{p} 3} =0.
\ee
We see that there are no stationary points in this case. This proves part (b) of Theorem \ref{Main2}. In this case, we can also describe completely the flow for diagonal initial data. The diagonal property is preserved, and setting $g_{\bar{b} a}(t)= \lambda_a(t) \delta_{ab}$, we find that $\lambda_3(t)=\lambda_3(0)$ is constant in time, while $\lambda_1(t)$ and $\lambda_2(t)$ satisfy the ODE system,
\be
\p_t \lambda_1 = {\sqrt{\lambda_3(0)} \over 2} \sqrt{ {\lambda_1 \over \lambda_2}}, \ \ \p_t \lambda_2 = {\sqrt{\lambda_3(0)} \over 2} \sqrt{ {\lambda_2 \over \lambda_1}}.
\ee
In particular,
\bea
\p_t{\lambda_1\over\lambda_2}=
0
\eea
Thus the ratio $\lambda_1(t)/\lambda_2(t)$ is constant for all time. Substituting in the previous equation, we can solve explicitly for $\lambda_1(t)$ and $\lambda_2(t)$,
\bea
\lambda_1(t)=\lambda_1(0)+{1\over 2}t\sqrt{\lambda_1(0)\lambda_3(0)\over\lambda_2(0)},
\qquad
\lambda_2(t)=\lambda_2(0)+{1\over 2}t\sqrt{\lambda_2(0)\lambda_3(0)\over\lambda_1(0)}
\eea
which tend both to $\infty$ as $t\to\infty$.

\subsubsection{Solvable case}

In this case, we may assume the Lie algebra satisfies the commutation relations 
\bea
[e_3,e_1]=e_1, \ [e_3,e_2]= -e_2, \ [e_1,e_2] = 0
\eea
that is, $c^1{}_{31}=1$, $c^2{}_{32} = -1$, and all the other structure constants vanish. Substituting these structure constants into the flow (\ref{flowmetric}) gives the following system
\bea \label{solvable-g11-flow}
\p_t g_{\bar{1} 1} &=& {1 \over 2 \| \Omega \|} \bigg( g^{3 \bar{3}} g_{\bar{1} 1} - \beta g^{3 \bar{3}}  g^{3 \bar{3}} g_{\bar{1} 1} \bigg),
\\
\label{solvable-g12-flow}
\p_t g_{\bar{1} 2} &=& {1 \over 2 \| \Omega \|} \bigg( g^{3 \bar{3}} g_{\bar{1} 2} - \beta g^{3 \bar{3}}  g^{3 \bar{3}} g_{\bar{1} 2} \bigg),
\\
\label{solvable-g13-flow}
\p_t g_{\bar{1} 3} &= &{1 \over 2 \| \Omega \|} \bigg( -g^{1 \bar{3}} g_{\bar{1} 1}  + g^{3 \bar{2}} g_{\bar{1} 2} + \beta g^{1 \bar{3}}  g^{3 \bar{3}} g_{\bar{1} 1} + \beta g^{3 \bar{2}} g^{3 \bar{3}} g_{\bar{1} 2}  \bigg),
\\
\p_t g_{\bar{2} 2} &=& {1 \over 2 \| \Omega \|} \bigg( g^{3 \bar{3}} g_{\bar{2} 2} - \beta g^{3 \bar{3}}  g^{3 \bar{3}} g_{\bar{2} 2}  \bigg)
\\
\label{solvable-g23-flow}
\p_t g_{\bar{2} 3} &=& {1 \over 2 \| \Omega \|} \bigg( g^{1 \bar{3}} g_{\bar{2} 1} - g^{2 \bar{3}} g_{\bar{2} 2} + \beta g^{1 \bar{3}}  g^{3 \bar{3}} g_{\bar{2} 1}  + \beta g^{2 \bar{3}}  g^{3 \bar{3}} g_{\bar{2} 2}  \bigg),
\\
\label{solvable-g33-flow}
\p_t g_{\bar{3} 3} &=& {1 \over 2 \| \Omega \|} \bigg( g^{1 \bar{1}} g_{\bar{1} 1} + g^{2 \bar{2}} g_{\bar{2} 2} - g^{2 \bar{1}} g_{\bar{1} 2} - g^{1 \bar{2}} g_{\bar{2} 1} \nonumber\\
&& - \beta g^{1 \bar{1}}  g^{3 \bar{3}} g_{\bar{1} 1}  - \beta g^{2 \bar{2}} g^{3 \bar{3}} g_{\bar{2} 2}  - \beta g^{2 \bar{1}} g^{3 \bar{3}} g_{\bar{1} 2} - \beta g^{1 \bar{2}} g^{3 \bar{3}} g_{\bar{2} 1} \bigg).
\eea
Here we set
\be \label{beta-defn}
\beta =  {\alpha' \tau \over 4}.
\ee
We begin by identifying the stationary metrics. 
Let $g_{\bar{a} b}$ be a stationary metric. From (\ref{solvable-g11-flow}) we see that $\beta>0$ and $g^{3 \bar{3}} = \beta^{-1}$. Substituting into (\ref{solvable-g13-flow}), we obtain
\be
2 g^{3 \bar{2}} g_{\bar{1} 2} = 0.
\ee
If $g_{\bar{1} 2} = 0$, the stationary point is of the desired form. Hence we must have $g^{3 \bar{2}}=0$. Similarly, substituting $g^{3 \bar{3}} = \beta^{-1}$ into (\ref{solvable-g23-flow}) leads to either $g_{\bar{2} 1} = 0$ or $g^{1 \bar{3}}=0$. Hence we may assume that $g^{3 \bar{2}}=g^{3 \bar{1}}=0$ which implies that $g_{\bar{3} 1}= g_{\bar{3} 2}=0$ and $g_{\bar{3} 3} = \beta$. Next, by (\ref{solvable-g33-flow}), we conclude
\be
0= 2(g^{2 \bar{1}} g_{\bar{1} 2} + g^{1 \bar{2}} g_{\bar{2} 1})= 4 \Re \, \{ g^{2 \bar{1}} g_{\bar{1} 2} \}.
\ee
By the formula for inverse matrices, we have
\be
g^{2 \bar{1}} g_{\bar{1} 2} = -{|g_{\bar{2} 1}|^2 g_{\bar{3} 3} \over \det g}.
\ee
Therefore $g_{\bar{1} 2} = 0$ and the solution is of the desired form.
\bigskip
\par It follows that the stationary metrics are exactly the metrics which satisfy
\be \label{stationary}
g_{\bar{1} 2}=0, \ \ g^{3 \bar{3}} = \beta^{-1}.
\ee
These equations can also be rewritten as
\bea
g_{\bar 12}=0,
\quad
{|g_{\bar 13}|^2\over
g_{\bar 11}}+{|g_{\bar 23}|^2\over g_{\bar 22}}
=g_{\bar 33}-\beta.
\eea
In particular, there are stationary points $g_{\bar{b} a}$ which are not diagonal. For example, setting $g_{\bar{1} 2}= g_{\bar{2} 1}=0$, $g_{\bar{3} 3} = 2 + \beta$ and all other entries to $1$ gives a stationary point. More generally, the moduli space of solutions requires locally two complex parameters $g_{\bar 13}$ and $g_{\bar 23}$, and two real parameters $g_{\bar 11}$, $g_{\bar 22}$.

\medskip
Next, we examine the Anomaly flow. First, we consider the case of initial metrics with $g_{\bar 12}(0)=0$. This condition is clearly preserved under the flow, and the flow for the other components of the metric becomes
\bea\label{g11}
\p_t g_{\bar{1} 1} &=& {1 \over 2 \| \Omega \|}  g^{3 \bar{3}} g_{\bar{1} 1}(1 - \beta g^{3 \bar{3}}) ,
\\
\label{g22}
\p_t g_{\bar{2} 2} &=& {1 \over 2 \| \Omega \|}  g^{3 \bar{3}} g_{\bar{2} 2}(1 - \beta g^{3 \bar{3}}) ,
\\
\label{g33}
\p_t g_{\bar{3} 3} &=& {1 \over 2 \| \Omega \|} (g^{1\bar{1}} g_{\bar{1} 1} (1-\beta g^{3 \bar{3}}) + g^{2 \bar{2}} g_{\bar{2} 2}(1-\beta g^{3 \bar{3}}) ),
\\
\label{g13}
\p_t g_{\bar{1} 3} &=& {1 \over 2 \| \Omega \|} g^{1 \bar{3}} g_{\bar{1} 1} (-1 + \beta g^{3 \bar{3}}) ,
\\
\label{g23}
\p_t g_{\bar{2} 3} &=& {1 \over 2 \| \Omega \|} g^{2 \bar{3}} g_{\bar{2} 2} (-1 + \beta g^{3 \bar{3}}).
\eea
We shall use the following simple formulas for the entries of the inverse metric $g^{a\bar b}$,
\bea \label{inverse-id}
&&
g^{1\bar 1} = {g_{\bar 22} g_{\bar 33} - |g_{\bar 23}|^2 \over \det g}, \ \ \ g^{2\bar 2} = {g_{\bar 11} g_{\bar 33} - |g_{\bar 13}|^2 \over \det g}, \ \ \ 
g^{3\bar 3} = {g_{\bar 11} g_{\bar 22}\over \det g}
\nonumber\\
&&
g^{1 \bar 3} = -{g_{\bar 22} g_{\bar 13} \over \det g}, \ \ \ g^{2\bar 3} = -{g_{\bar 11} g_{\bar 23} \over \det g}.
\eea
We note the following identities:

\begin{itemize}
 \item $g_{\bar 22} = a g_{\bar 11}$ with $a= {g_{\bar 22}\over g_{\bar 11}}(0)$. This follows directly from equation (\ref{g11}) and (\ref{g22}).
 \item $|g_{\bar 13}| = b g_{\bar 11}$ with $b= {|g_{\bar 13}|\over g_{\bar 11}}(0)$. Indeed, putting $g^{1\bar 3}$ into equation (\ref{g13}) and $g^{3 \bar{3}}$ into (\ref{g11}), we obtain
 \bea
\p_t g_{\bar{1} 3} = {1 \over 2 \| \Omega \|} {g_{\bar 13} g_{\bar 22}g_{\bar 11}\over \det g} (1 - \beta g^{3 \bar{3}}), \ \ \p_t g_{\bar{1} 1} = {1 \over 2 \| \Omega \|}  {g_{\bar 11} g_{\bar 22}g_{\bar{1} 1}\over \det g}(1 - \beta g^{3 \bar{3}}).
 \eea
Hence
\bea
\p_t \ln | g_{\bar 13}| = \p_t \ln g_{\bar 11}
\eea
and this implies the desired relation.
\item $|g_{\bar 23}| = c g_{\bar 11}$ with $c= {|g_{\bar 23}|\over g_{\bar 11}}(0)$. The proof is similar to the previous case.
\item $g^{3\bar 3} = {d\over g^2_{\bar 11}}$ with $d= g^2_{\bar 11}(0)g^{3\bar 3}(0)$. First, we compute
\bea \label{p-t-g33}
\p_t g^{3 \bar{3}} &=& {g_{\bar{2} 2} \p_t g_{\bar{1} 1} \over \det g} + {g_{\bar{1} 1} \p_t g_{\bar{2} 2} \over \det g} - {g_{\bar{1} 1} g_{\bar{2} 2} \over (\det g)^2} \p_t \det g \nonumber\\
&=& {1 \over \| \Omega \|} g^{3 \bar{3}} g^{3 \bar{3}} (1 - \beta g^{3 \bar{3}}) - {g_{\bar{1} 1} g_{\bar{2} 2} \over (\det g)^2} \p_t \det g.
\eea
It follows from (\ref{g13}) and (\ref{g23}) that 
\be
\p_t |g_{\bar{1} 3}|^2 = {1 \over \|\Omega \|} g^{3 \bar{3}} ( 1 - \beta g^{3 \bar{3}}) |g_{\bar{1} 3}|^2, \ \ \p_t |g_{\bar{2} 3}|^2 = {1 \over \|\Omega \|} g^{3 \bar{3}} ( 1 - \beta g^{3 \bar{3}}) |g_{\bar{2} 3}|^2.
\ee
Next, 
\bea
\p_t \, \det g &=& \p_t ( g_{\bar{1} 1} g_{\bar{2} 2} g_{\bar{3} 3} - g_{\bar{1} 1} |g_{\bar{2} 3}|^2 - g_{\bar{2} 2} | g_{\bar{1} 3}|^2) \nonumber \\
&=& \p_t g_{\bar{1} 1} g_{\bar{2} 2} g_{\bar{3} 3} + g_{\bar{1} 1} \p_t g_{\bar{2} 2} g_{\bar{3} 3} + g_{\bar{1} 1} g_{\bar{2} 2} \p_t g_{\bar{3} 3} - \p_t g_{\bar{1} 1} |g_{\bar{2} 3}|^2 \nonumber\\
&& g_{\bar{1} 1} \p_t |g_{\bar{2} 3}|^2 - \p_t g_{\bar{2} 2} | g_{\bar{1} 3}|^2 - g_{\bar{2} 2} \p_t | g_{\bar{1} 3}|^2 \nonumber\\
&=& {1 \over 2 \| \Omega \|} ( 1 - \beta g^{3 \bar{3}}) \left( 2 g^{3 \bar{3}} g_{\bar{1} 1} g_{\bar{2} 2} g_{\bar{3} 3} + g_{\bar{1} 1} g_{\bar{2} 2} (g^{1 \bar{1}} g_{\bar{1} 1} + g^{2 \bar{2}} g_{\bar{2} 2}) \right) \nonumber\\
&& - {3 \over 2 \| \Omega \|} g^{3 \bar{3}} ( 1 - \beta g^{3 \bar{3}}) ( g_{\bar{1} 1} |g_{\bar{2} 3}|^2 + g_{\bar{2} 2} |g_{\bar{1} 3}|^2). 
\eea
Using (\ref{inverse-id}) yields
\bea \label{p-t-detg}
\p_t \, \det g &=& {1 \over 2 \| \Omega \|} ( 1 - \beta g^{3 \bar{3}}) g^{3 \bar{3}} \bigg( 2 g_{\bar{1} 1} g_{\bar{2} 2} g_{\bar{3} 3} + (\det g) (g^{1 \bar{1}} g_{\bar{1} 1} + g^{2 \bar{2}} g_{\bar{2} 2}) \nonumber\\
&& -3 g_{\bar{1} 1} |g_{\bar{2} 3}|^2 -3 g_{\bar{2} 2} |g_{\bar{1} 3}|^2 \bigg) \nonumber\\
&=& {2 \over  \| \Omega \|} ( 1 - \beta g^{3 \bar{3}}) g^{3 \bar{3}} \det g.
\eea
Combining (\ref{p-t-g33}) and (\ref{p-t-detg}),
\bea
\p_t g^{3\bar 3} = - {1\over \|\Omega\|} g^{3\bar 3} g^{3\bar 3} (1- \beta g^{3\bar 3}).
\eea
Comparing this equation with (\ref{g11}), it follows that
\bea
\p_t g^{3\bar 3} = -2 {\p_t g_{\bar 11} \over g_{\bar 11} }g^{3\bar 3}, \ \ \p_t \ln g^{3\bar 3} = 2 \p_t \ln {1\over g_{\bar 11}},
\eea
which implies $g^{3\bar 3} = {d\over g^2_{\bar 11}}$ with $d= g^2_{\bar 11}(0)g^{3\bar 3}(0)$.
\end{itemize}

Denote $\lambda= g_{\bar 11}$. Using the previously derived formulas and solving for $g_{\bar{3} 3}$ from the expression for $\det g$, it follows that
\be \label{metric-lambda}
g_{\bar{1} 1} = \lambda, \ \ g_{\bar{2} 2} = a \lambda, \ \ |g_{\bar{1} 3}| = b \lambda, \ \ |g_{\bar{2} 3}| = c \lambda,  \ \ \det g = {a \over d} \, \lambda^4, \ \ g_{\bar{3} 3} = {\lambda^2 \over d} + \bigg( {c^2 \over a} + b \bigg) \lambda.
\ee
Putting the above relations into equation (\ref{g11}), 
\bea
\p_t \lambda = {(\det g)^{1/2}\over 2} g^{3\bar 3} g_{\bar 11} (1- \beta g^{3\bar 3})={1\over 2} \sqrt{a\over d} \lambda^2 \cdot {d\over \lambda^2} \cdot \lambda (1- {\beta d \over \lambda^2})
\eea
Then
\bea
\p_t \lambda^2 = \sqrt{ad}\,  (\lambda^2 - \beta d),
\eea
and hence
\bea
|\lambda^2(t) - \beta d| = |\lambda^2(0) - \beta d| \cdot e^{\sqrt{ad}\,t}.
\eea
Suppose $\beta>0$. Then $\lambda$ is uniformly bounded above and away from zero as $t \rightarrow - \infty$, hence from (\ref{metric-lambda}) we see that the metric $g_{\bar{b} a}$ is uniformly bounded and remains non-degenerate. Thus the flow exists for all time $t<0$ and $g^2_{\bar 11}(t) \rightarrow \beta d$ as $t\rightarrow -\infty$ for any initial data. In particular, this implies that $g^{3\bar 3} \rightarrow \beta^{-1}$. This completes the description of the flow for initial data satisfying the condition $g_{\bar 12}(0)=0$.

\bigskip
We consider now the case of initial data with $g_{\bar{1} 2}\not=0$.
In this case, we claim that the flow cannot converge to a non-degenerate metric.
Indeed, if $\p_t g_{\bar{b} a} \rightarrow 0$ then $g_{\bar{1} 2} \rightarrow 0$. However, from (\ref{solvable-g11-flow}) and (\ref{solvable-g12-flow}), we deduce
\be
{\p_t g_{\bar{1} 1} \over g_{\bar{1} 1}} = {\p_t |g_{\bar{1} 2}| \over |g_{\bar{1} 2}|}.
\ee
Therefore for all times, 
\be
|g_{\bar{1} 2}| = C g_{\bar{1} 1}, \ \ C = {|g_{\bar{1} 2}(0)| \over g_{\bar{1} 1}(0)}.
\ee
The convergence to $0$ of $g_{\bar 12}$ implies then the convergence to $0$
of $g_{\bar{1} 1} \rightarrow 0$, which contradicts 
the requirement that the limit be a non-degenerate metric.
In particular, given a stationary solution $g_\infty$, if we perturb it in the $g_{\bar{1} 2}$ direction and run the Anomaly flow, the flow will not take the metric back to $g_\infty$. This completes the proof of part (c) of Theorem \ref{Main2}.

\subsubsection{The Semi-Simple Case}

This is the case of $X=SL(2,{\bf C})$, whose standard basis $\{ e_a \}$ has structure constants $c^k{}_{ij} = \epsilon_{kij}$, the Levi-Civita symbol. 

We begin by showing that the metrics $g_{\bar ab}=2\beta\delta_{ab}$ are the only stationary points of the flow, where $\beta$ is defined in (\ref{beta-defn}). In particular, the existence of a stationary point requires $\beta>0$.

Assume then that $g_{\bar{a} b}$ is a stationary point of the flow. Let $P$ be a matrix such that $\bar{P}^{\bar{a}}{}_{\bar{p}} g_{\bar{a} b} P^b{}_{q} = \delta_{\bar{p} q}$, and perform a change of basis by defining $f_i = e_r P^r{}_{i}$, $[f_a,f_b] = f_c k^c{}_{ab}$. As previously discussed (\ref{metric-vs-P}), in this new frame $\{ f_a \}$ we have $\omega = i f^a \wedge \bar{f}^a$. The fixed point equation for the flow (\ref{flowmetric}) becomes
\be \label{fixedpt-kappa-eqn}
h_{\bar{a} b} = \beta \, \overline{k^m{}_{ap}} h_{\bar{m} s} k^s{}_{bp},
\ee
where we defined
\be
h_{\bar{a} b} = \overline{k^m{}_{ap}} k^m{}_{b p}.
\ee
If we let $\kappa_{\bar{a} b} =  \overline{c^m{}_{ap}} c^m{}_{b p}$, by symmetry of the Levi-Civita symbol we obtain $\kappa_{\bar{a} b} = 2 \delta_{ab}$. By using the transformation laws (\ref{transformation}), we derive
\be \label{killingform-sl2c}
k^\ell{}_{ij} = (P^{-1})^\ell{}_s P^r{}_i P^q{}_j \epsilon_{srq}, \ \ h_{\bar{a} b} = \bar{P}^{\bar{q}}{}_{\bar{a}} P^p{}_{b} \kappa_{\bar{q} p} = 2 \bar{P}^{\bar{r}}{}_{\bar{a}} P^r{}_{b}.
\ee
If we substitute (\ref{killingform-sl2c}) into (\ref{fixedpt-kappa-eqn}) and use $g^{a \bar{b}} = P^a{}_{r} \bar{P}^{\bar{b}}{}_{\bar{r}}$, cancellation occurs and we are left with
\be\label{killingform-sl2c1}
h_{\bar{a} b} = 2 \beta \, \bar{P}^{\bar{p}}{}_{\bar{a}} P^q{}_b g^{r \bar{s}} \epsilon_{\ell p s} \epsilon_{\ell q r}.
\ee 
We note the formula
\be\label{levi-civita}
\sum_{\ell}\epsilon_{\ell ps} \epsilon_{\ell qr} = \delta_{pq}\delta_{sr} - \delta_{pr} \delta_{qs}.
\ee
Combining (\ref{killingform-sl2c1}) and (\ref{levi-civita}), we obtain
\bea
h_{\bar a b} &=& 2\beta \, \bar{P}^{\bar{p}}{}_{\bar{a}} P^q{}_b g^{r \bar{s}} \,( \delta_{pq} \delta_{sr} - \delta_{pr} \delta_{qs})\nonumber
\\
&=&
2\beta  \, (\bar{P}^{\bar{r}}{}_{\bar{a}} P^r{}_{b}) {\rm Tr}(g^{-1})- 2\beta\, \bar{P}^{\bar{p}}{}_{\bar{a}} g^{p \bar{q}} P^q{}_b.
\eea
Substituting the relations $h_{\bar{a} b} = 2 \bar{P}^{\bar{r}}{}_{\bar{a}} P^r{}_{b}$ and $g^{a \bar{b}} = P^a{}_{r} \bar{P}^{\bar{b}}{}_{\bar{r}}$, it follows that
\bea
h_{\bar a b} &=& \beta\, h_{\bar ab} (P^m{}_r \bar P^{\bar m}{}_{\bar r})- 2\beta\, \bar{P}^{\bar p}{}_{\bar a} P^p{}_{m} \bar{P}^{\bar q}{}_{\bar m} P^q{}_b
= {1\over 2} \beta\, h_{\bar ab} {\rm Tr}\,h- {1\over 2} \beta\, h_{\bar a m} h_{\bar m b}.
\eea
Hence
\bea \label{kappaequation}
h_{\bar a b} (\beta {\rm Tr} h -2) = \beta h_{\bar a m} h_{\bar m b}.
\ea
By multiplying $h^{b\bar a}$ on both sides of the equation, we obtain $\beta\, h_{\bar b b} = \beta {\rm Tr} h-2$. In particular, ${\rm Tr} h = 3\beta^{-1}$. Putting this back into equation (\ref{kappaequation}),
\bea
h_{\bar a b} = \beta h_{\bar a m}h_{\bar m b}.
\eea
Thus, $\beta h_{\bar ab}$ is an invertible idempotent matrix. This implies that $\beta\, h=I$ and hence $\bar P^{\bar r}{}_{\bar a} P^r{}_{b}={1\over 2} h_{\bar ab} = {1\over 2 \beta} \delta_{ab}$. Therefore,
\bea
g_{\bar{a} b} = (\bar{P}^{-1})^{\bar{r}}{}_{\bar{a}} (P^{-1})^r{}_{b} = 2 \beta\, \delta_{ab}
\eea
as was to be shown.

\medskip
Next, we show that the stationary metric is asymptotically unstable
\footnote{Recall that a stationary point for a flow is said to be asymptotically stable if the flow will converge to the stationary point for any initial data in some neighborhood of the point. The flow is said to a asymptotically unstable if it is not asymptotically stable.}. For this, it suffices to show that the flow can be restricted to a submanifold of Hermitian metrics, and that restricted to this submanifold, the flow is asymptotically unstable. We choose this submanifold to be the submanifold of metrics $g_{\bar ab}$ which are diagonal with respect to the given basis of invariant holomorphic vector fields,
\bea
g_{\bar ab}=\lambda_b \delta_{ab}.
\eea
Using the explicit form (\ref{flowmetric}) of the flow and the fact that the structure constants are given by $\e_{abc}$, it is easy to verify that the diagonal form of metrics is preserved along the flow. The flow reduces to the following ODE system for the eigenvalues $\lambda_a$,
\bea
\p_t \lambda_1 &=& { (\lambda_1 \lambda_2 \lambda_3)^{1/2} \over 2} \bigg( \lambda_2^{-1} \lambda_3 +  \lambda_3^{-1} \lambda_2  - 2 \beta \, \lambda_1^{-1} - \beta \,\lambda_1 (\lambda_2)^{-2}   - \beta \, \lambda_1 (\lambda_3)^{-2}  \bigg),
\nonumber
\\
\p_t \lambda_2 &=& { (\lambda_1 \lambda_2 \lambda_3)^{1/2} \over 2} \bigg( \lambda_1^{-1} \lambda_3 +  \lambda_3^{-1} \lambda_1  - 2 \beta \, \lambda_2^{-1} - \beta \,\lambda_2 (\lambda_1)^{-2}   - \beta \, \lambda_2 (\lambda_3)^{-2}  \bigg),
\nonumber
\\
\p_t \lambda_3 &=& { (\lambda_1 \lambda_2 \lambda_3)^{1/2}  \over 2} \bigg( \lambda_1^{-1} \lambda_2 +  \lambda_2^{-1} \lambda_1  - 2 \beta \, \lambda_3^{-1} - \beta \,\lambda_3 (\lambda_1)^{-2}   - \beta \, \lambda_3 (\lambda_2)^{-2}  \bigg).
\nonumber
\eea

The linearization of the flow at the stationary point $\lambda_1=\lambda_2=\lambda_3=2\beta$ is easily worked out,
\bea
\p_t\lambda_a=\sum_b Q_{ab}(\delta\lambda_b)
\eea
where the matrix $Q=(Q_{ab})$ is given by
\bea
Q=\sqrt{\beta\over 2}\,\left[ {\begin{array}{ccc} 
0 &1 &1 \\
1 & 0 &1 \\
1 & 1 & 0 \\
\end{array}}\right]
\eea
The matrix $Q$ has eigenvalues $-\sqrt{\beta/2}$ and $2\sqrt{\beta/2}$ with multiplicities $2$ and $1$ respectively. By a classical theorem on ordinary differential equations (see e.g. Theorem 3.3 in \cite{V}),
the presence of an eigenvalue with strictly positive real part implies that the flow is asymptotically unstable. Part (d)
of Theorem \ref{Main2} has now been proved, completing the proof of the theorem.

\bigskip

\subsubsection{Remarks}

We conclude with several remarks.

\medskip

$\bullet$ For simplicity, we have formulated Theorem \ref{Main2} under the assumption that $\alpha'\tau>0$. The behavior of the Anomaly flow can be readily worked out as well 
by similar methods when $\alpha'\tau\leq 0$. The arguments are in fact simpler in that case, because there is then no cancellation between the two terms on the right hand side (\ref{flowmetric}) of the flow. We leave the details to the reader.

\medskip
$\bullet$ In general, the sign of the right hand side in the Anomaly flow is dictated by the requirement that the flow be weakly parabolic. But in the case of Lie groups, the flow reduces to an ODE, and both signs are allowed. The opposite sign can be obtained from the sign we chose here simply by a time-reversal.

\medskip

$\bullet$ The remaining remarks are about the semi-simple case $SL(2,{\bf C})$. The eigenvalues of the linerarized operator at the stationary point imply that there is a stable surface and an unstable curve near this point. The stable surface appears difficult to identify explicitly, but the unstable curve is easily found. It is given by the line of metrics proportional to the identity matrix, $g_{\bar ab}=\lambda\,\delta_{ab}$. This line is preserved under the flow, which reduces to
\bea
\p_t\lambda=\lambda^{3\over 2} \left( 1-{2\beta\over\lambda} \right).
\eea 
This equation can be solved explicitly by
\be
\lambda(t) = 2 \beta \, \bigg( { C e^{\sqrt{2\beta} t} + 1 \over 1 - C e^{\sqrt{2\beta} t}} \bigg)^2, \ \ {\rm if} \, \lambda(0) > 2 \beta,
\ee
\be
\lambda(t) = 2 \beta \, \bigg( { 1 - C e^{\sqrt{2\beta} t} \over C e^{\sqrt{2\beta} t} + 1 } \bigg)^2, \ \ {\rm if} \, \lambda(0) < 2 \beta,
\ee
where 
\bea
C = \bigg| {\sqrt{\lambda}(0) - \sqrt{2 \beta} \over \sqrt{\lambda}(0) + \sqrt{2 \beta} } \bigg| < 1.
\eea
This shows that the flow terminates in finite time at $T={1\over \sqrt{2\beta}}\log {1\over C}$, with $\lambda(t)\to+\infty$
as $t\to T$ if $\lambda(0)>2\beta$, and $\lambda(t)\to 0$ as $t\to T$
if $\lambda(0)<2\beta$.

\medskip

$\bullet$ More generally, if two eigenvalues are equal at some time, then they are equal for all time. This follows from rewriting the flow as
\bea
\p_t\lambda_1
=
-{(\lambda_1\lambda_2\lambda_3)^{1\over 2}\over 2}
\bigg(\beta({2\over\lambda_1}+{\lambda_1\over\lambda_2^2}+{\lambda_1\over\lambda_3^2})-{\lambda_2^2+\lambda_3^2\over\lambda_2\lambda_3}\bigg)
\eea
with similar formulas for $\p_t\lambda_2$ and $\p_t\lambda_3$. In particular, we have
\bea
\label{pt}
\p_t\log {\lambda_1\over\lambda_2}
&=&
-{(\lambda_1\lambda_2\lambda_3)^{1\over 2}\over 2}
\big(\beta({1\over \lambda_1^2}-{1\over\lambda_2^2})
-
{\lambda_2^2-\lambda_1^2\over \lambda_1\lambda_2\lambda_3}\big)
\nonumber\\
&=&
-{(\lambda_1\lambda_2\lambda_3)^{1\over 2}\over 2}
{\lambda_2^2-\lambda_1^2\over\lambda_1^2\lambda_2^2}
\big(\beta-{\lambda_1\lambda_2\over \lambda_3}\big).
\eea
Let $[0,T)$ be the maximum time of existence of the flow. This equation implies that if any two eigenvalues are equal at some time $t_0$, then they are identically equal on the whole interval $[0,T)$. Indeed, by the Cauchy-Kowalevska theorem, the eigenvalues are analytic functions of $t$ near any time where they are all strictly positive. The equation (\ref{pt}) implies that, if say $\lambda_1$ and $\lambda_2$ are equal at $t_0$, then all derivatives in time of $\lambda_1$ and $\lambda_2$ at $t_0$ are also equal, as we can see by differentiating the equation (\ref{pt}). By analyticity, they must be equal in a neighborhood of $t_0$. Thus the set where $\lambda_1$ and $\lambda_2$ coincide is both open and closed. This establishes our claim. 

\medskip
$\bullet$ It follows that if the eigenvalues at the initial time are ordered as
\bea
\lambda_1\geq\lambda_2\geq\lambda_3
\eea
then this ordering is preserved by the flow. The configuration space can be divided into the invariant and mutually disjoint subsets $\{\lambda_1>\lambda_2>\lambda_3\}$, $\{\lambda_1=\lambda_2>\lambda_3\}$, $\{\lambda_1>\lambda_2=\lambda_3\}$, $\{\lambda_1=\lambda_2=\lambda_3\}$. We have already shown that the flow diverges to $+\infty$ on the last invariant subset. We shall next make a few remarks on the other sets.

\medskip
$\bullet$ On each of the other invariant subsets, we have the following: 
(a) If $\lambda_1(0)<2\beta$, then $\lambda_1(t)$ is monotone decreasing, and in particular less than $\lambda_1(0)$ for all time $t\in [0,T)$;
(b) If $\lambda_3(0)>2\beta$, then $\lambda_3(t)$ is monotone increasing, and in particular greater than $\lambda_3(0)$ for all time $t\in [0,T)$.
To see this, we express the flow as
\bea
\p_t\lambda_1
=
{(\lambda_1\lambda_2\lambda_3)^{1\over 2}\over 2}
\bigg
({\lambda_2\over\lambda_3}
+
{\lambda_3\over\lambda_2}-\beta({2\over\lambda_1}+{\lambda_1\over\lambda_2^2}
+
{\lambda_1\over\lambda_3^2})\bigg).
\eea
We shall make use of the following two estimates
\bea
\label{upper}
{\lambda_1\over\lambda_2}
+
{\lambda_2\over\lambda_1}-\beta({2\over\lambda_3}+{\lambda_3\over\lambda_1^2}
+
{\lambda_3\over\lambda_2^2})
\geq
{\lambda_1\over\lambda_2}
+
{\lambda_2\over\lambda_1}-\beta({2\over\lambda_3}+{\lambda_3\over\lambda_3^2}
+
{\lambda_3\over\lambda_3^2})
=
{\lambda_1\over\lambda_2}
+
{\lambda_2\over\lambda_1}-4\beta{1\over\lambda_3}
\eea
and
\bea
\label{lower}
{\lambda_2\over\lambda_3}
+
{\lambda_3\over\lambda_2}-\beta({2\over\lambda_1}+{\lambda_1\over\lambda_2^2}
+
{\lambda_1\over\lambda_3^2})
\leq
{\lambda_2\over\lambda_3}
+
{\lambda_3\over\lambda_2}-\beta {2\over\lambda_2}-\beta{2\over \lambda_3}
=
{1\over\lambda_2}(\lambda_3-2\beta)
+
{1\over\lambda_3}(\lambda_2-2\beta).
\nonumber\\
\eea
We can now establish (a). First, we claim that $\lambda_1(t)<2\beta$ for any time $t\in [0,T)$. Otherwise, let $t_0$ be the first time when $\lambda_1(t_0)=2\beta$. 
Then $\lambda_1(t)<2\beta$ on the interval $[0,t_0)$. 
On the interval $[0,t_0]$, we have then $\lambda_3\leq\lambda_2\leq \lambda_1\leq 2\beta$, and the inequality (\ref{lower}) implies that $\p_t\lambda_1<0$ on this interval. It follows that $\lambda_1(t_0)\leq \lambda_1(0)<2\beta$, which contradicts our assumption. But now that we know that $\lambda_1(t)<2\beta$ for all time $t$, the same inequality (\ref{lower}) shows that $\lambda_1(t)$ is a strictly monotone decreasing function of time.

Next, we establish (b). Again, 
let $t_0$ be the first time when $\lambda_3(t_0)=2\beta$. On the interval $[0,t_0)$, we can apply the inequality (\ref{upper}) and obtain
\bea
\p_t\lambda_3
\geq 
{(\lambda_1\lambda_2\lambda_3)^{1\over 2}\over 2}
\bigg({\lambda_1\over\lambda_2}+{\lambda_2\over\lambda_1}
-\beta{4\over\lambda_3}\bigg)
>
0
\eea
where we used the inequality ${\lambda_1\over\lambda_2}+{\lambda_2\over\lambda_1}\geq 2$. It follows that $\lambda_3(t_0)>\lambda_3(0)>2\beta$, which is a contradiction. Thus $\lambda_3(t)$ is a strictly monotone increasing function of time.

\bigskip
\noindent
{\bf Acknowledgements} The authors would especially like to thank Teng Fei for many stimulating discussions.

\bigskip
Department of Mathematics, Columbia University, New York, NY 10027, USA

\smallskip

phong@math.columbia.edu

\bigskip
Department of Mathematics, Columbia University, New York, NY 10027, USA

\smallskip
picard@math.columbia.edu

\bigskip
Department of Mathematics, University of California, Irvine, CA 92697, USA

\smallskip
xiangwen@math.uci.edu

\end{document}